\documentclass[12pt]{article}

\usepackage{amsmath, amsfonts, amssymb, amsthm, subfigure, geometry}

\usepackage[dvips, all, knot]{xy}
\xyoption{all}

\usepackage{fancyhdr}
\usepackage{sectsty}
\allsectionsfont{\normalsize\centering}

\newcommand{\SL}{{\rm SL}}
\newcommand{\GL}{{\rm GL}}
\newcommand{\Sp}{{\rm Sp}}
\newcommand{\GSp}{{\rm GSp}}

\newcommand{\Z}{{\mathbb Z}}
\newcommand{\F}{{\mathbb F}}

\newcommand{\Card}{{\rm Card}}
\newcommand{\diag}{{\rm diag}}
\newcommand{\lk}{{\rm lk}}
\newcommand{\ord}{{\rm ord}}
\newcommand{\ti}{{\rm totally isotropic}}

\newcommand{\cB}{{\cal B}}
\newcommand{\cO}{{\cal O}}

\newtheorem{nlem}{Lemma}[section]
\newtheorem{nprop}{Proposition}[section]
\newtheorem{nthm}{Theorem}[section]
\newtheorem{ncor}{Corollary}[section]

\theoremstyle{definition}
\newtheorem*{rem}{Remark}

\CompileMatrices

\begin{document}
\title{Distance in the Affine Buildings of $\SL_n$ and $\Sp_n$}
\author{Alison Setyadi \\
\normalsize{setyadi@member.ams.org}}
\date{}
\maketitle
\thispagestyle{empty}
\begin{abstract}
For a local field $K$ and $n \geq 2$, let $\Xi_n$ and $\Delta_n$ denote
the affine buildings naturally associated to the special linear and
symplectic groups $\SL_n(K)$ and $\Sp_n(K)$, respectively. We relate the
number of vertices in $\Xi_n$ ($n \geq 3$) close (i.e., gallery distance
$1$) to a given vertex in $\Xi_n$ to the number of chambers in $\Xi_n$
containing the given vertex, proving a conjecture of Schwartz and
Shemanske. We then consider the special vertices in $\Delta_n$ ($n \geq
2$) close to a given special vertex in $\Delta_n$ (all the vertices in
$\Xi_n$ are special) and establish analogues of our results for
$\Delta_n$.
\end{abstract}

\section*{Introduction}

\begin{figure}
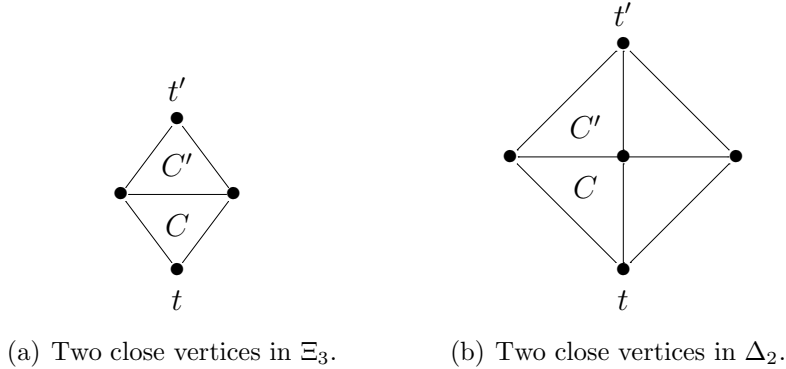

\centering
\subfigure[Two close vertices in $\Xi_3$.]{
\label{SL3_close}
\xy 
(0, 0)*{\bullet}="c4";
(0, -4)*{\text{$t$}};
(0, 6)*{\text{$C$}};
(-7.5, 10)*{\bullet}="u3";
(-22.5, 10)*{};
(7.5, 10)*{\bullet}="u4";
(22.5, 10)*{};
(0, 20)*{\bullet}="uu4";
(0, 24)*{\text{$t'$}};
(0, 14)*{\text{$C'$}};
(0, -9)*{};
{\ar@{-} "c4";"u3" };
{\ar@{-} "c4";"u4" };
{\ar@{-} "u3";"u4" };
{\ar@{-} "uu4";"u3" };
{\ar@{-} "uu4";"u4" };
\endxy
}
\hspace{0.4in}
\subfigure[Two close vertices in $\Delta_2$.]{
\label{Sp2_close}
\xy 
(0, 30)*{\bullet}="t";
(0, 34)*{\text{$t'$}};
(-5, 19)*{\text{$C'$}};
(-15, 15)*{\bullet}="cl";
(-22.5, 15)*{};
(0, 15)*{\bullet}="c";
(15, 15)*{\bullet}="cr";
(22.5, 15)*{};
(0, 0)*{\bullet}="b";
(0, -4)*{\text{$t$}};
(-5, 11)*{\text{$C$}};
(0, -9)*{};
{\ar@{-} "t";"c" };
{\ar@{-} "t";"cr" };
{\ar@{-} "t";"cl" };
{\ar@{-} "cl";"cr" };
{\ar@{-} "b";"c" };
{\ar@{-} "b";"cr" };
{\ar@{-} "b";"cl" };
\endxy}
\caption[Examples of close vertices]{Examples of close vertices.}
\end{figure}

A building is a finite-dimensional simplicial complex in which any two
of its chambers (maximal simplices) can be connected by a gallery. In
other words, if $\Delta$ is a building, then for any chambers $C, D \in
\Delta$, there is a sequence $C = C_0, C_1, \ldots, C_m = D$ of chambers
in $\Delta$ such that $C_i$ and $C_{i + 1}$ are adjacent (share a
codimension-one face) for all $0 \leq i \leq m - 1$; in this case, the
number $m$ is the length of the gallery $C_0, \ldots, C_m$. The
combinatorial distance between $C$ and $D$
is the minimal length of a gallery in $\Delta$ connecting $C$ and $D$
(see \cite[p.\ 14]{brown}). Following \cite[p.\ 15]{brown},
define the {\it distance} between any non-empty simplices $A, B \in
\Delta$ to be the minimal length of a gallery in $\Delta$ starting at a
chamber containing $A$ and ending at a chamber containing $B$ (cf.\ 
\cite[p.\ 125]{ss}). Then the vertices $t, t' \in \Delta$ are distance
one apart or {\it close} if and only if there are adjacent chambers $C,
C' \in \Delta$ such that $t \in C$, $t' \in C'$, but $t, t' \not\in C
\cap C'$ (the simplex shared by $C$ and $C'$); i.e., if and only if $t$
and $t'$ are in adjacent chambers in $\Delta$ but not a common one
(cf.\ \cite[p.\ 127]{ss}). Figures \ref{SL3_close} and \ref{Sp2_close}
show close vertices in the affine buildings naturally associated to
$\SL_3(K)$ and $\Sp_2(K)$, respectively, for any local field $K$. Note
that if $\Delta$ is a building and $t, t' \in \Delta$ are close
vertices, then as vertices in the underlying graph of $\Delta$, $t$ and
$t'$ are not graph distance $1$ apart but are always graph distance $2$
apart.

Let $K$ be a local field with valuation ring $\cO$, uniformizer $\pi$,
and residue field $k \cong \F_q$, and let $\Xi_n$ denote the affine building
naturally associated to $\SL_n(K)$. In \cite[Theorem 3.3]{ss}, Schwartz
and Shemanske show that for all $n \geq 3$, the number $\omega_n$ of
vertices in $\Xi_n$ close to a given vertex in $\Xi_n$ is the number of
right cosets of $\GL_n(\cO)$ in $\GL_n(\cO)\diag(1, \pi, \ldots, \pi,
\pi^2)\GL_n(\cO)$; i.e., the Hecke operator $\GL_n(\cO)\diag(1, \pi,
\ldots, \pi, \pi^2)\GL_n(\cO)$ acts as a generalized adjacency operator
on $\Xi_n$. They also conjecture that for all $n \geq 3$, $q \cdot r_n =
r_{n - 2}\ \omega_n$, where $r_n$ is the number of chambers in $\Xi_n$
containing a given vertex, with $r_1 := 1$ (see the remark following
\cite[Proposition 3.4]{ss}).

In Section \ref{SLn_close}, we prove Schwartz and Shemanske's conjecture
in two ways. Our first approach is via module theory. More precisely, we
use the description of the chambers in $\Xi_n$ in terms of lattices in
an $n$-dimensional $K$-vector space (see, for example, \cite[p.\ 
115]{ronan}) to obtain an explicit formula for $\omega_n$ (Proposition
\ref{SLnd1count}); together with Schwartz and Shemanske's formula for
$r_n$ \cite[Proposition 2.4]{ss}, this proves Theorem \ref{SLnrelation}.
Our second approach is through combinatorics (Theorem
\ref{SLnrelation2}). Specifically, we show that if $t, t' \in \Xi_n$ are
close vertices, then there is a one-to-one correspondence between the
galleries of length $1$ in $\Xi_n$ whose initial chamber contains $t$
and whose ending chamber contains $t'$ and the chambers in the spherical
$A_{n - 3}(k)$ building. This gives an explanation for the relationship
between $\omega_n$ and $r_n$ in terms of the structure of $\Xi_n$. In
Section \ref{Spn_close}, we consider the special vertices in the affine
building $\Delta_n$ naturally associated to $\Sp_n(K)$ ($n \geq 2$)
close to a given special vertex in $\Delta_n$ (all the vertices in
$\Xi_n$ are special). Using the fact that $\Delta_n$ is a subcomplex of
$\Xi_{2n}$, we adapt the proofs of the results for close vertices in
$\Xi_{2n}$ to prove analogues for $\Delta_n$. In particular, we
establish analogues of \cite[Theorem 3.3]{ss} and Theorem
\ref{SLnrelation} (Theorems \ref{3.3} and \ref{Spnrelation},
respectively) and a partial analogue of Theorem \ref{SLnrelation2}
(Proposition \ref{Spnrelation2}). Note that while every vertex in
$\Xi_{2n}$ is special, only two vertices in each chamber in $\Delta_n$
are special; hence, our analysis for $\Delta_n$ requires more care than
that needed for $\Xi_{2n}$.

%
%
After proving Theorems \ref{SLnrelation} and \ref{Spnrelation},
we learned that
the formulas in Propositions \ref{SLnd1count} and \ref{Spnd1count} are
both
special cases of a result of Parkinson \cite[Theorem 5.15]{parkinsonja}
and that the formula in Proposition \ref{SLnd1count} also follows from a
result of Cartwright \cite[Lemma 2.2]{cartwright}. We view the buildings
$\Xi_n$ and $\Delta_n$ as combinatorial objects naturally associated to
$\SL_n(K)$ and $\Sp_n(K)$, respectively, and make use of the lattice
descriptions of these buildings (see \cite{garrett} and \cite{ronan}).
As a result, our methods require little more than
the definition of a building---namely, some module theory. In contrast
to our approach, Cartwright views $\Xi_n$ in terms of hyperplanes,
affine transformations, and convex hulls,
and Parkinson considers
buildings via
root systems and Poincar\'{e} polynomials of Weyl groups. The
numbers $\omega_n$ and $\omega(\Delta_n)$ that we use are special cases
of Parkinson's
$N_\lambda$, which he uses
to define vertex set averaging operators on arbitrary locally finite,
regular
affine buildings
and whose
formula he uses
to prove results about those
operators.

I thank Paul Garrett for the idea behind the proof of Proposition
\ref{SLnd1count}, and hence that of Proposition \ref{Spnd1count}.
Finally, the results contained here form part of my doctoral thesis,
which I wrote under the guidance of Thomas R.\ Shemanske.

\section{Close Vertices in the Affine Building $\Xi_n$ of $\SL_n(K)$}
\label{SLn_close}

From now on, $K$ is a local field with discrete valuation ``ord,''
valuation ring $\cO$, uniformizer $\pi$, and residue field $k \cong
\F_q$. For any finite-dimensional $K$-vector space $V$, define a
{\it lattice} in $V$ to be a free $\cO$-submodule of $V$ of rank
$\dim_K V$, with two lattices $L$ and $L'$ in $V$ {\it homothetic} if
$L' = \alpha L$ for some $\alpha \in K^\times$; write $[L]$ for the
homothety class of the lattice $L$.


The
affine building $\Xi_n$ naturally associated to
$\SL_n(K)$ can be modeled as an $(n - 1)$-dimensional simplicial complex
as follows
(see \cite[p.\ 115]{ronan}). Let
$V$ be an $n$-dimensional $K$-vector space. Then a {\it vertex}
in $\Xi_n$ is a homothety class of lattices in $V$, and two
vertices $t, t' \in \Xi_n$ are {\it incident}
if there are representatives $L \in t$ and $L' \in t'$ such that $\pi
L \subseteq L' \subseteq L$; i.e., such that $L'/\pi L$ is a
$k$-subspace of $L/\pi L$. Thus,
a
{\it chamber} (maximal simplex)
in $\Xi_n$
has $n$ vertices $t_0, \ldots, t_{n - 1}$ with
representatives $L_i \in t_i$ such that $\pi L_0 \subsetneq L_1
\subsetneq \cdots \subsetneq L_{n - 1} \subsetneq L_0$ and
$[L_1 : \pi L_0] = q = [L_i : L_{i - 1}]$ for all $2 \leq i \leq n -
1$.
From now on,
write that a chamber in $\Xi_n$ corresponds to the chain
$\pi L_0 \subsetneq L_1 \subsetneq \cdots \subsetneq L_{n - 1}
\subsetneq L_0$
only when the lattices $L_0, \ldots, L_{n - 1}$ satisfy the conditions
in the last sentence.

\begin{figure}
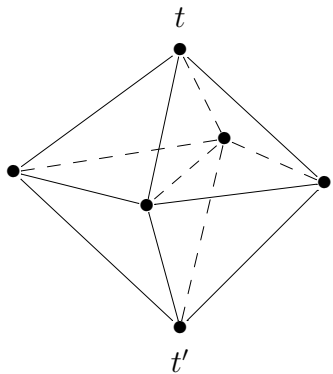

\[
\xy 0;/r0.35pc/:
(0, 0)*{\bullet}="top";
(0, 3)*{\text{$t$}};
(-3, -14)*{\bullet}="front";
(4, -8)*{\bullet}="back";
(0, -25)*{\bullet}="bottom";
(0, -28)*{\text{$t'$}};
(-15, -11)*{\bullet}="left";
(13, -12)*{\bullet}="right";
{\ar@{-} "top";"front" };
{\ar@{--} "top";"back" };
{\ar@{--} "front";"back" };
{\ar@{-} "front";"bottom" };
{\ar@{--} "back";"bottom" };
{\ar@{-} "left";"front" };
{\ar@{-} "left";"bottom" };
{\ar@{-} "left";"top" };
{\ar@{--} "left";"back" };
{\ar@{-} "right";"top" };
{\ar@{-} "right";"bottom" };
{\ar@{-} "right";"front" };
{\ar@{--} "right";"back" };
\endxy
\]
\caption{Two close vertices in $\Xi_4$.}
\label{closeSL4}
\end{figure}

For the rest of this section, $n \geq 3$. Let
$t \in \Xi_n$ be
a vertex with representative $L$. Then a chamber $C \in \Xi_n$
containing $t$ corresponds to a chain of
the form
\begin{equation} \label{SLnchamberchain}
\pi L \stackrel{q}{\subsetneq} L_1 \stackrel{q}{\subsetneq} \cdots
\stackrel{q}{\subsetneq} L_{n - 1} \stackrel{q}{\subsetneq} L
\end{equation}
(cf.\ 
\cite[p.\ 323]{garrett}).
The codimension-one face in $C$ not containing $t$
thus corresponds to the chain
\[
L_1 \stackrel{q}{\subsetneq} \cdots \stackrel{q}{\subsetneq} L_{n - 1},
\]
and a vertex in $\Xi_n$ is close to $t$ if it has a representative $M
\neq L$ such that
\begin{equation} \label{SLnclosechain}
\pi M \stackrel{q}{\subsetneq} L_1 \stackrel{q}{\subsetneq} \cdots
\stackrel{q}{\subsetneq} L_{n - 1} \stackrel{q}{\subsetneq} M.
\end{equation}
Given the lattices $L_1$ and $L_{n - 1}$, the possible $L$ and $M$
satisfy $L_{n - 1} \subsetneq L \neq M \subsetneq \pi^{-1}L_1$. On the
other hand, if $t, t' \in \Xi_n$ are close vertices, then there must be
representatives $L \in t$ and $M \in t'$ and lattices $L_1, \ldots, L_{n
- 1}$ as in \eqref{SLnchamberchain} 
such that $L_{n - 1} \subsetneq L \neq M \subsetneq \pi^{-1}L_1$. Recall
that if $M_1$ and $M_2$ are free, rank $n$, $\cO$-modules with $M_1
\subseteq M_2$, then $M_1 \subseteq M' \subseteq M_2$ implies $M'$ is
also a free, rank $n$, $\cO$-module. Thus, both
$L \cap M$ and $L + M$ are lattices in $V$. Furthermore, $L \neq M$
and $[L : L_{n - 1}] = q = [M : L_{n - 1}]$
imply
$L \cap M = L_{n - 1}$ and $L + M = \pi^{-1}L_1$, but we can vary
$L_2, \ldots, L_{n - 2}$ as long as $L_1 \subsetneq L_2 \subsetneq
\cdots \subsetneq L_{n - 2} \subsetneq L_{n - 1}$. In other words, if
$t$ and $t'$ are close vertices in $\Xi_n$, there may be two (or more)
pairs of adjacent chambers $C$ and $C'$ in $\Xi_n$ with $t \in C$, $t'
\in C'$, but $t, t' \not\in C \cap C'$ (see Figure \ref{closeSL4}).
We return to this later.

Before we count the number of vertices in $\Xi_n$ close to a given
vertex $t \in \Xi_n$, we make a few observations. Fix a representative
$L \in t$. Since
$L/\pi L \cong k^n$,
the Correspondence Theorem and the fact that any
$\cO$-submodule
of $L$ containing $\pi L$
is a lattice in $V$ imply that
the number of $L_1$ is
the number of $1$-dimensional $k$-subspaces of $L/\pi L$. Similarly,
given $L_1$ as above, the number of lattices $L_{n - 1}$ with
$L_1 \subsetneq L_{n - 1} \subsetneq L$ and $[L : L_{n - 1}] = q$ is
the number of $(n - 2)$-dimensional $k$-subspaces of $L/L_1 \cong k^{n
- 1}$. Finally, given $L_1$ and $L_{n - 1}$ as above, the number of
lattices $M \neq L$ such that $L_{n - 1} \subsetneq M \subsetneq
\pi^{-1}L_1$ is one less than the number of non-trivial, proper
$k$-subspaces of $\pi^{-1}L_1/L_{n - 1} \cong k^2$.
\begin{nprop} \label{SLnd1count}
If $t \in \Xi_n$ is a vertex, then the number $\omega_n$ of vertices
in $\Xi_n$ close to $t$ is
\[
\frac{q^n - 1}{q - 1} \cdot \frac{q^{n - 1} - 1}{q - 1} \cdot q
\]
(independent of $t$).
\begin{proof}
This follows from the preceding
comments,
duality,
and the fact that the number of $1$-dimensional subspaces of $\F_q^m$
is exactly $(q^m - 1)/(q - 1)$.
\end{proof}
\end{nprop}
\begin{ncor}
The number of right cosets of $\GL_n(\cO)$ in $\GL_n(\cO)\diag(1, \pi,
\ldots, \pi, \pi^2)\GL_n(\cO)$ is
$((q^n - 1)(q^{n - 1} - 1) \cdot q)/(q - 1)^2$.
\begin{proof}
This follows from \cite[Theorem 3.3]{ss} and the last proposition.
\end{proof}
\end{ncor}

Let $r_n$ be the number of chambers in $\Xi_n$ containing a vertex $t
\in \Xi_n$. Then
\cite[Proposition 2.4]{ss} and
the last proposition establish
the conjecture following Proposition 3.4 of \cite{ss}:
\begin{nthm} \label{SLnrelation}
For all $n \geq 3$, $q \cdot r_n = r_{n - 2}\ \omega_n$, where $r_1 =
1$.
\end{nthm}

We now use the structure of $\Xi_n$ to give a combinatorial
proof for the relationship given in Theorem \ref{SLnrelation}. Fix a
vertex $t \in \Xi_n$. Then we can try to
count the number of vertices in $\Xi_n$ close to $t$ by counting the
number of galleries (in $\Xi_n$) of length $1$ starting at a chamber
containing $t$ and ending at a chamber not containing $t$. By
definition,
there are $r_n$ chambers $C \in \Xi_n$ containing $t$. Since a
chamber in $\Xi_n$
adjacent to $C$ and not containing $t$ must contain the
codimension-one face
in $C$ not containing $t$,
\cite[p.\ 324]{garrett} implies that
there are $q$
chambers in $\Xi_n$
adjacent to $C$ not containing $t$; hence,
there are exactly $r_n \cdot q$ galleries of length $1$ in $\Xi_n$
whose initial chamber contains $t$ and whose ending chamber does not
contain $t$. On the other hand, if
$t' \in \Xi_n$ is a vertex close to $t$,
we count $t'$ more than once if there is
more than one gallery of length $1$ in $\Xi_n$ whose initial chamber
contains $t$ and whose ending chamber contains $t'$ (see Figure
\ref{closeSL4}); hence,
$\omega_n = (r_n \cdot q)/m(t, t')$, where $m(t, t')$ is the number of
galleries of length $1$ in $\Xi_n$ whose initial chamber contains $t$
and whose ending chamber contains $t'$.

To determine
$m(t, t')$,
fix the following notation for the rest of this section. For close
vertices
$t, t' \in \Xi_n$,
let $L \in t$, $M \in t'$ be representatives such that there are
lattices $L_1, \ldots, L_{n - 1}$ as in \eqref{SLnchamberchain} and
\eqref{SLnclosechain}. 
Recall that
$L_1 = \pi(L + M)$ and $L_{n - 1} = L \cap M$, but we can vary $L_2,
\ldots, L_{n - 2}$ as long as $L_1 \subsetneq L_2 \subsetneq \cdots
\subsetneq L_{n - 2} \subsetneq L_{n - 1}$. Since any gallery $C, C'$
in $\Xi_n$ such that $C = \{t, [L_1], \ldots, [L_{n - 1}]\}$ and $C' =
\{t', [L_1], \ldots, [L_{n - 1}]\}$ satisfies $C \cap C' = \{[L_1],
\ldots, [L_{n - 1}]\}$, each
gallery in $\Xi_n$ counted by $m(t, t')$ is uniquely determined by the
lattices $L_2, \ldots, L_{n - 2}$.
Define two vertices in $\Xi_n$ to be
{\it adjacent} if they are distinct and incident.
\begin{nprop} \label{adjacentreps}
Let $t, t' \in \Xi_n$ be adjacent vertices. If $L \in t$, then there
is a unique representative $L' \in t'$ such that $\pi L \subsetneq L'
\subsetneq L$.
\begin{proof}
Since
$t$ and $t'$ are incident and $t \neq t'$, there are representatives
$M \in t$ and $M' \in t'$ such that
$\pi M \subsetneq M' \subsetneq M$. Moreover, $M$ and $L$ are
homothetic, so
$L = \alpha M$ for some $\alpha \in K^\times$; hence, $\pi L
\subsetneq \alpha M' \subsetneq
L$. Let $L' = \alpha M'$. If
$L'' \in t'$ such that $\pi L \subsetneq L'' \subsetneq L$, let
$\beta \in K^\times$ such that $L'' = \beta L'$. Suppose $\ord(\beta) = m$.
Then $\pi L \subsetneq L' \subsetneq L$ implies
$\pi^{m + 1} L \subsetneq L'' \subsetneq \pi^m L$ and
$L = \pi^m L$; i.e.,
$L''
= L'$.
\end{proof}
\end{nprop}

Consider the set of vertices in $\Xi_n$ that are adjacent to $t, t',
[L + M]$, and $[L \cap M]$
(in the case $n = 3$, this set is empty), and define
two such vertices to be incident
if they are incident
as vertices in $\Xi_n$.
Let $\Xi_n^c(t, t')$ be the set consisting of
\begin{itemize}
\item the empty set,
\item all vertices in $\Xi_n$ adjacent to $t, t', [L + M]$, and $[L
\cap M]$, and
\item all finite sets $A$ of vertices in $\Xi_n$ adjacent to $t, t',
[L + M]$, and $[L \cap M]$ such that any two vertices in $A$ are
adjacent.
\end{itemize}
Then $\Xi_n^c(t, t')$ is a simplicial complex. In particular,
$\Xi_n^c(t, t')$
is a subcomplex of $\Xi_n$.

\begin{nlem} \label{SLnclosesimplex}
If $\emptyset \neq A \in \Xi_n^c(t, t')$ is an $i$-simplex, then $A$
corresponds to a chain of lattices $M_1 \subsetneq \cdots \subsetneq
M_{i + 1}$, where $\pi(L + M) \subsetneq M_1 \subsetneq \cdots
\subsetneq M_{i + 1} \subsetneq L \cap M$. In particular, $A$ has at
most $n - 3$ vertices.
\begin{proof}
We proceed by induction on $i$. If $i = 0$, then $A$
adjacent to $[L \cap M]$ implies $A$
has a unique representative $M_1$
such that $\pi(L \cap M) \subsetneq M_1
\subsetneq L \cap M$ by Proposition \ref{adjacentreps}. Then by
\cite[p.\ 322]{garrett}, either $M_1
\subsetneq \pi(L + M)$ or $M_1
\supsetneq \pi(L + M)$. In the second case, we are done, so assume
$M_1
\subsetneq \pi(L + M)$. Then $\pi(L \cap M) \subsetneq M_1
\subsetneq \pi(L + M)$. On the other hand, $\pi(L \cap M) \subsetneq
\pi L \subsetneq \pi(L + M)$ and $[\pi(L + M) : \pi(L \cap M)] =
q^2$. Since $A$
is adjacent to $t$, \cite[p.\ 322]{garrett} implies that either $M_1
\subsetneq \pi L$ or $M_1
\supsetneq \pi L$. Thus, either $\pi(L \cap M) \subsetneq M_1
\subsetneq \pi L \subsetneq \pi(L + M)$ or $\pi(L \cap M) \subsetneq
\pi L \subsetneq M_1 \subsetneq \pi(L + M)$, which is impossible given
the previous index computation.
%

Now suppose $0 \leq i \leq n - 5$ and that the claim holds for any
$i$-simplex in $\Xi_n^c(t, t')$. Let $A \in \Xi_n^c(t, t')$ be an $(i
+ 1)$-simplex
and
$x \in A$
a vertex.
Then
the $i$-simplex $A - \{x\}$ corresponds to a chain of lattices $M_1'
\subsetneq \cdots \subsetneq M_{i + 1}'$ such that $\pi(L + M)
\subsetneq M_1' \subsetneq \cdots \subsetneq M_{i + 1}' \subsetneq L
\cap M$.
By the last paragraph,
$x$ has a representative $M'$ such that $\pi(L + M) \subsetneq M'
\subsetneq L \cap M$. If $M' \subsetneq M_1'$, set $M_1 = M'$ and $M_j
= M_{j - 1}'$ for all $2 \leq j \leq i + 2$. Otherwise, $M' \supsetneq
M_1'$ by \cite[p.\ 322]{garrett}. Let $j \in \{1, \ldots, i + 1\}$ be
maximal such that $M' \supsetneq M_j'$.
If $j = i + 1$, set $M_\ell = M_\ell'$ for all $1 \leq \ell \leq i +
1$ and $M_{i + 2} = M'$. Setting
$M_\ell = M_\ell'$ for all $1 \leq \ell \leq j$, $M_{j + 1} =
M'$, and $M_\ell = M_{\ell - 1}'$ for all $j + 2 \leq \ell \leq i + 2$
finishes the proof if $j \neq i + 1$. Finally, note that if the claim
holds for $i \geq n - 3$, then $A$ corresponds to a chain of lattices
$M_1 \subsetneq \cdots \subsetneq M_{i + 1}$, where $\pi(L + M)
\subsetneq M_1 \subsetneq \cdots \subsetneq M_{i + 1} \subsetneq L
\cap M$, contradicting the fact that $[L \cap M : \pi(L + M)] = q^{n -
2}$. 
\end{proof}
\end{nlem}

Write $\Xi_n^s(k)$ for the spherical $A_n(k)$ building described in
\cite[p.\ 4]{ronan}.
\begin{nprop}
For any close vertices $t, t' \in \Xi_n$,
$\Xi_n^c(t, t')$
is isomorphic (as a poset) to
$\Xi_{n - 3}^s(k)$ (independent of $t$ and $t'$), where
$\Xi_0^s(k) = \emptyset$.
\begin{proof}
Let $L \in t, M \in t'$ be as in the paragraph preceding Proposition
\ref{adjacentreps}, and
let $\Xi_{n - 3}^s(k)$ be the spherical $A_{n - 3}(k)$ building with
simplices
the empty set, together with the nested sequences
of non-trivial, proper $k$-subspaces of $(L \cap M)/\pi(L + M)$. Then
by the Correspondence Theorem and
the last lemma,
there is a bijection between the $i$-simplices in $\Xi_n^c(t, t')$ and
the $i$-simplices in $\Xi_{n - 3}^s(k)$ for all $i$. Since
this bijection preserves the partial order (face) relation, it
is a poset isomorphism.
\end{proof}
\end{nprop}
\begin{nthm} \label{SLnrelation2}
If $t, t' \in \Xi_n$ are close vertices, then
$m(t, t') =
r_{n - 2}$ (independent of $t$ and $t'$). In particular,
$\omega_n = (r_n \cdot q)/r_{n - 2}$.
\begin{proof}
By the
last proposition and
previous comments, $m(t, t')$ is
the number of chambers in
$\Xi_{n - 3}^s(k)$. The proof now follows from
\cite[Proposition 2.4]{ss}.
\end{proof}
\end{nthm}
\section{Close Vertices in the Affine Building $\Delta_n$ of
$\Sp_n(K)$} \label{Spn_close}

Let $\Delta_n$ denote the affine building naturally associated to
$\Sp_n(K)$. Then $\Delta_n$ is a subcomplex of $\Xi_{2n}$, and there
is a natural embedding of $\Delta_n$ in $\Xi_{2n}$. As we will see,
this embedding allows us to derive information about $\Delta_n$ and to
prove results for $\Delta_n$ by adapting the proofs of the analogous
results for $\Xi_{2n}$. As noted in the introduction, while all the
vertices in $\Xi_{2n}$ are special, only two vertices in each chamber
in $\Delta_n$ are special. Consequently, the $\Sp_n$ case requires
more care than that needed in the last section. We start by looking at
properties of $\Delta_n$ that we need to consider close vertices in
$\Delta_n$.

\subsection{The building $\Delta_n$}

The building
$\Delta_n$ can be modeled as an $n$-dimensional simplicial complex as
follows (see
\cite[pp.\ 336 -- 337]{garrett}). Fix a $2n$-dimensional $K$-vector
space $V$ endowed with a non-degenerate, alternating bilinear form
$\langle \cdot, \cdot \rangle$, and recall that a subspace $U$ of $V$
is {\it totally isotropic} if $\langle u, u' \rangle = 0$ for all $u,
u' \in U$. A lattice $L$ in $V$ is {\it primitive}
if $\langle L, L \rangle
\subseteq \cO$ and $\langle \cdot, \cdot \rangle$ induces a
non-degenerate,
alternating $k$-bilinear form
on $L/\pi L$. Then a {\it vertex}
in $\Delta_n$ is a
homothety class
of lattices in $V$ with a representative $L$ such that
there is a primitive lattice $L_0$ with $\langle L, L \rangle
\subseteq \pi\cO$ and $\pi L_0 \subseteq L \subseteq L_0$;
equivalently, $L/\pi L_0$ is a \ti\ $k$-subspace of $L_0/\pi L_0$. Two
vertices $t, t' \in \Delta_n$ are {\it incident}
if there are representatives $L \in t$
and $L' \in t'$ such that there is
a primitive lattice $L_0$ with $\langle L, L \rangle \subseteq
\pi\cO$, $\langle L', L' \rangle \subseteq \pi\cO$, and
either $\pi L_0 \subseteq L \subseteq L' \subseteq L_0$ or $\pi L_0
\subseteq L' \subseteq L \subseteq L_0$. Thus,
a
{\it chamber}
in $\Delta_n$ has $n + 1$ vertices $t_0, \ldots, t_n$ with
representatives $L_i \in t_i$ such that $L_0$ is primitive, $\langle
L_i, L_i \rangle \subseteq \pi\cO$
for all $1 \leq i \leq n$,
and
$\pi L_0 \subsetneq L_1 \subsetneq \cdots \subsetneq L_n \subsetneq
L_0$. 
From now on,
write that a chamber in $\Delta_n$ corresponds to the chain
$\pi L_0 \subsetneq L_1 \subsetneq \cdots \subsetneq L_n \subsetneq
L_0$
only when the lattices $L_0, \ldots, L_n$ satisfy the conditions in
the last sentence.

Recall that a basis $\{u_1, \ldots, u_n, w_1, \ldots, w_n\}$ for $V$
is {\it symplectic} if $\langle u_i, w_j \rangle = \delta_{ij}$
(Kronecker delta) and $\langle u_i, u_j \rangle = 0 = \langle w_i, w_j
\rangle$ for all $i, j$.
If a
$2$-dimensional, \ti\ subspace $U$ of $V$ is a
hyperbolic plane,
then a frame
is an unordered $n$-tuple $\{\lambda_1^1, \lambda_1^2\}, \ldots,
\{\lambda_n^1, \lambda_n^2\}$ of pairs of lines ($1$-dimensional
$K$-subspaces) in $V$ such that
\begin{enumerate}
\item $\lambda_i^1 + \lambda_i^2$ is a hyperbolic plane for all $1
\leq i \leq n$,
\item $\lambda_i^1 + \lambda_i^2$ is orthogonal to $\lambda_j^1 +
\lambda_j^2$ for all $i \neq j$, and
\item $V = (\lambda_1^1 + \lambda_1^2) + \cdots + (\lambda_n^1 +
\lambda_n^2)$.
\end{enumerate}
A vertex $t \in \Delta_n$ lies in the apartment specified by the frame
$\{\lambda_1^1, \lambda_1^2\}, \ldots, \{\lambda_n^1, \lambda_n^2\}$
if for any representative $L \in t$, there are lattices $M_i^j$ in
$\lambda_i^j$ for all $i, j$ such that $L = M_1^1 + M_1^2 + \cdots +
M_n^1 + M_n^2$. The following lemma is easily established.
\begin{nlem} \label{symplecticframe}
$\mbox{}$
\begin{enumerate}
\item Every symplectic basis for $V$ specifies an apartment of
$\Delta_n$.
\item If $\Sigma$ is an apartment of $\Delta_n$, there is a symplectic
basis $\{u_1, \ldots, u_n, w_1, \ldots, w_n\}$ for $V$ such that every
vertex in $\Sigma$ has the form
\[
[\cO\pi^{a_1}u_1 + \cdots + \cO\pi^{a_n}u_n + \cO\pi^{b_1}w_1 + \cdots
+ \cO\pi^{b_n}w_n]
\]
for some $a_i, b_i \in \Z$.
\end{enumerate}
\end{nlem}
\begin{rem}
A frame specifying an apartment of $\Delta_n$ also specifies an
apartment of $\Xi_{2n}$ (see \cite[p.\ 323]{garrett}). In particular,
a symplectic basis for $V$ specifies an apartment of $\Xi_{2n}$.
\end{rem}       

Since
$\pi$ is fixed, if $\cB = \{u_1, \ldots, u_n, w_1, \ldots, w_n\}$ is a
symplectic basis for $V$, follow \cite[p.\ 3411]{shemanske}
and write $(a_1, \ldots, a_n; b_1, \ldots, b_n)_{\cB}$
for the lattice $\cO\pi^{a_1}u_1 + \cdots + \cO\pi^{a_n}u_n +
\cO\pi^{b_1}w_1 + \cdots + \cO\pi^{b_n}w_n$ and $[a_1, \ldots, a_n;
b_1, \ldots, b_n]_{\cB}$
for its homothety class.
%
Then the lattice $L = (a_1, \ldots, a_n; b_1, \ldots, b_n)_{\cB}$ is
primitive if and only if $a_i + b_i = 0$ for all $i$ by \cite[p.\ 
3411]{shemanske}, and $[L]$ is a {\it special} vertex in $\Delta_n$ if
and only if $a_i + b_i = \mu$ is constant for all $i$ by
\cite[Corollary 3.4]{shemanske}. Note that by \cite[p.\ 3412]{shemanske},
a chamber in $\Delta_n$ has exactly two special vertices.
\begin{nlem} \label{primitiverep}
Let $t \in \Delta_n$ be a vertex with a primitive representative $L$,
and let $\Sigma$ be an apartment of $\Delta_n$ containing $t$. Then
there is a symplectic basis $\cB$ for $V$ specifying $\Sigma$ as in
Lemma \ref{symplecticframe} such that $L = (0, \ldots, 0; 0, \ldots,
0)_{\cB}$.
\begin{proof}
This follows from 
Lemma \ref{symplecticframe} and
\cite[p.\ 3411]{shemanske}.
\end{proof}
\end{nlem}

Let $t \in \Delta_n$ be a vertex. Then
the link of $t$ in $\Delta_n$, 
denoted $\lk_{\Delta_n} t$, is a building (see \cite[Proposition
IV.1.3]{brown}) that
is isomorphic (as a poset) to the subposet of $\Delta_n$
consisting of those simplices containing $t$ by \cite[p.\ 31]{brown}.
In particular, if $A \in \Delta_n$ is a codimension-one simplex
containing $t$ and $A' \in \lk_{\Delta_n} t$ is the codimension-one
simplex corresponding to $A$, then the number of chambers in
$\Delta_n$ containing $A$ is the number of chambers in $\lk_{\Delta_n}
t$ containing $A'$.
Note that if $t$ is {\it special}, then \cite[p.\ 35]{tits} implies
$\lk_{\Delta_n} t$ is isomorphic to the spherical $C_n(k)$ building
$\Delta_n^s(k)$ described in \cite[pp.\ 5 -- 6]{ronan}.
\begin{nprop} \label{chambercount}
Every special vertex in $\Delta_n$ is contained in exactly
$r(\Delta_n) = \prod_{m = 1}^n\ ((q^{2m} - 1)/(q - 1))$
chambers in $\Delta_n$.
\begin{proof}
Let $t \in \Delta_n$ be a special vertex. By the preceding comments
and \cite[pp.\ 5 -- 6]{ronan}, it suffices to count the number of
maximal flags of non-trivial, \ti\ subspaces of a $2n$-dimensional
$k$-vector space endowed with a non-degenerate, alternating bilinear
form. An obvious modification of the proof of \cite[Proposition
2.4]{ss} finishes the proof.
\end{proof}
\end{nprop}
\begin{rem}
The number $r(\Delta_n)$ in the last proposition corresponds to the
number $r_n$ given in \cite[Proposition 2.4]{ss}.
Since $\Sp_1(K) = \SL_2(K)$,
set $r(\Delta_1) = q + 1$ for completeness.
\end{rem}
\begin{nprop} \label{codimno}
If $A \in \Delta_n$ is a codimension-one simplex, then $A$ is
contained in
exactly $q + 1$ chambers in $\Delta_n$.
\begin{proof}
Let $t$ be a special vertex in $A$ and
$A'$
the codimension-one simplex in $\lk_{\Delta_n} t$ corresponding to
$A$. By the comments preceding the last proposition, it suffices to
count
the number of chambers in
$\Delta_n^s(k)$ containing $A'$. A case-by-case analysis finishes the
proof.
\end{proof}
\end{nprop}

We now use the fact that $\Delta_n$ is a subcomplex of $\Xi_{2n}$ to
derive information about $\Delta_n$.
For a vertex $t \in \Xi_{2n}$ with representative $L = \cO v_1 +
\cdots + \cO v_{2n}$ and $g \in \GL_{2n}(K)$, define $gt
= [\cO (gv_1) + \cdots + \cO (gv_{2n})]$. Then
$\GL_{2n}(K)$ acts transitively on the lattices in $V$.

Let
\[
J_n = \left(\begin{smallmatrix}
0 & I_n \\
-I_n & 0
\end{smallmatrix}\right) \text{ and }
\GSp_n(K) = \left\{g \in M_{2n}(K) : g^t J_n g = \nu(g)J_n
\textrm{ for some $\nu(g) \in K^\times$} \right\},
\]
so that
$\Sp_n(K)$
consists of the matrices $g \in \GSp_n(K)$ with $\nu(g) = 1$.
%
Alternatively, abuse notation and think of $\GSp_n(K)$ as
\[
\{g \in \GL_K(V) : \text{$\forall\ v_1, v_2 \in V,\ \exists\ \nu(g)
\in K^\times$ such that $\langle gv_1, gv_2 \rangle = \nu(g)\langle
v_1, v_2 \rangle$}\}.
\]
If $g \in \GL_{2n}(K)$ and
$\cB = \{v_1, \ldots, v_{2n}\}$ is a basis for $V$, write $g\cB$
for $\{gv_1, \ldots, gv_{2n}\}$.
\begin{nlem}
\label{Spprim}
The group $\Sp_n(K)$ acts on the set of primitive lattices in $V$.
\begin{proof}
Let $L$ be a primitive lattice in $V$, and let
$\Sigma$ be an apartment of $\Delta_n$ containing $[L]$ and $\cB$ a
symplectic basis for $V$ specifying $\Sigma$ as in Lemma
\ref{symplecticframe}. Then
$L = (a_1, \ldots, a_n; -a_1, \ldots, -a_n)_{\cB}$ by \cite[p.\ 
3411]{shemanske}; hence, for
$g \in \Sp_n(K)$,
$g\cB$
a symplectic basis for $V$
implies that $gL$ is primitive.
\end{proof}
\end{nlem}

For the rest of this section, let $\cB_0 = \{e_1, \ldots, e_n, f_1,
\ldots, f_n\}$ be the standard symplectic basis for $V$ ($f_i = e_{n +
i}$ for all $i$),
$L_0 = (0, \ldots, 0; 0, \ldots, 0)_{\cB_0}$, and $t_0 = [L_0]$.
Following \cite[p.\ 116]{ronan}, assign {\it types} to the vertices in
$\Xi_{2n}$ as follows: assign type $0$ to $t_0$ and type $\ord(\det g)
\mod 2n$ to
any other vertex $t
= [L] \in \Xi_{2n}$,
where $g \in \GL_{2n}(K)$ such that $L = gL_0$.
%
This induces
a labelling on the vertices in $\Delta_n$.
%
For the rest of this section, let
$C_0$ be the chamber in $\Delta_n$ whose vertices are the homothety
classes of the lattices
\begin{equation} \label{fundreps}
L_0 = (0, \ldots, 0; 0, \ldots, 0)_{\cB_0}, L_1 = (0, 1, \ldots, 1; 1,
\ldots, 1)_{\cB_0}, \ldots, L_n = (0, \ldots, 0; 1, \ldots,
1)_{\cB_0}.
\end{equation}
Note that
$[L_i]$ has type $2n - i$ for all $1 \leq i \leq n$. Recall that since
$\Delta_n$ is the affine building naturally associated to $\Sp_n(K)$,
$\Sp_n(K)$ acts on the vertices in $\Delta_n$ in a type-preserving
manner and also acts transitively on the chambers in $\Delta_n$.
\begin{nprop} \label{Spntypes}
If $t \in \Delta_n$ is a vertex, then $t$ has type $i$ for some $i
\equiv n, \ldots, 2n \mod 2n$.
\begin{proof}
By the preceding comments, it
suffices to show that for all $0 \leq j \leq n$, $[L_j]$ (as in
\eqref{fundreps}) has type $i$ for some $i \equiv n, \ldots, 2n \mod
2n$, which we already observed.
\end{proof}
\end{nprop}

We now use types to characterize the vertices in $\Delta_n$ with a
primitive representative, as well as those that are special.
\begin{nprop} \label{primitivetype}
A vertex in $\Delta_n$ has a primitive representative if and only if
it has type $0$.
\begin{proof}
Let $t \in \Delta_n$ be a type $0$ vertex and $C \in \Delta_n$ a
chamber containing $t$. Choose
$g \in \Sp_n(K)$ such that $gC_0 = C$. Then
$gL_0 \in t$. Since
$L_0$ is primitive,
Lemma \ref{Spprim} implies that $gL_0$ is primitive. Conversely, let
$t \in \Delta_n$ be
a vertex with a primitive representative $L$, and let $C \in \Delta_n$
be a chamber containing $t$. Let
$g \in \Sp_n(K)$ such that $gC = C_0$. Then
$gL = \pi^mL_j$ for some $0 \leq j \leq n$ and some $m \in \Z$.
If $L_j = (a_1, \ldots, a_n; b_1, \ldots, b_n)_{\cB_0}$ as in
\eqref{fundreps}, then
$gL =
(a_1 + m, \ldots, a_n + m; b_1 + m, \ldots, b_n + m)_{\cB_0}$. But
$gL$
primitive (by Lemma
\ref{Spprim})
implies that
$a_i + b_i = -2m$ for all $i$. By
\eqref{fundreps},
$m = 0$
and $gt = [L_0]$; hence,
$t$ has type $0$.
\end{proof}
\end{nprop}
%
\begin{nprop} \label{specialtype}
A vertex in $\Delta_n$ is special if and only if it has type $0$ or
$n$.
\begin{proof}
Let $t \in \Delta_n$ be a type $0$ (resp., type $n$) vertex, and let
$C \in \Delta_n$ be a chamber containing $t$. If
$g \in \Sp_n(K)$ such that $gC_0 = C$, then
$t = g[L_0]$ (resp., $t = g[L_n]$),
and $t$ is special by \cite[Corollary 3.4]{shemanske}.
%
Conversely, let $t \in \Delta_n$ be a special vertex. Let $C \in
\Delta_n$ be a chamber
containing $t$, $\Sigma$ an apartment of $\Delta_n$ containing $C$,
and $\cB$ a symplectic basis for $V$ specifying $\Sigma$ as in Lemma
\ref{symplecticframe}. By \cite[Corollary 3.4]{shemanske},
$t = [a_1, \ldots, a_n; \mu - a_1, \ldots, \mu - a_n]_{\cB}$ for some
$\mu \in \Z$. If
$g \in \Sp_n(K)$ such that $gC = C_0$, then
$gt =
[L_i]$ for some $0 \leq i \leq n$; hence,
$gt$ special, \cite[Corollary 3.4]{shemanske},
and \eqref{fundreps} imply
$i = 0$ or $i = n$, and
$t$ has type $0$ or $n$.
\end{proof}
\end{nprop}

We now consider the action of $\GSp_n(K)$ on the vertices in
$\Xi_{2n}$.
\begin{nprop} \label{action_type}
If $[L]$ is a type $i$ vertex in $\Xi_{2n}$, then for any $g \in
\GL_{2n}(K)$, the vertex $g[L] \in \Xi_{2n}$ has type $i + \ord(\det
g) \mod 2n$.
\begin{proof}
Since $[L]$ has type $i$, we can write $L = g_iL_0$, where $g_i \in
\GL_{2n}(K)$ with $\ord(\det g_i) \equiv i \mod 2n$. Then
$g[L]$
has type $\ord(\det(gg_i)) \mod 2n \equiv
i + \ord(\det g) \mod 2n$.
\end{proof}
\end{nprop}
\begin{ncor} \label{nonspecial}
If $g \in \GSp_n(K)$ with $\ord(\nu(g)) \equiv 1 \mod 2$, then $g$
maps a non-special vertex in $\Delta_n$ to a vertex in $\Xi_{2n}$
that is not in $\Delta_n$.
\begin{proof}
First note that $g \in \GSp_n(K)$ with $\ord(\nu(g)) \equiv 1 \mod 2$
implies $\ord(\det g) \equiv n \mod 2n$. If
$t$ is
a non-special vertex in $\Delta_n$, then
$t$ has type $i$ for some $n + 1 \leq i \leq 2n - 1$ by Propositions
\ref{Spntypes} and \ref{specialtype}. Thus, the last proposition
implies
$gt$ has type $i + n \mod 2n \in \{1, \ldots, n - 1\}$. Proposition
\ref{Spntypes} finishes the proof.
\end{proof}
\end{ncor}

\subsection{The building $\Delta_n$ in the building $\Xi_{2n}$}

Let $C \in \Delta_n$ be a chamber
corresponding to the chain
$\pi L_0 \subsetneq L_1 \subsetneq \cdots \subsetneq L_n \subsetneq
L_0$. Let $\Sigma$ be an apartment of $\Delta_n$ containing $C$, $\cB$
a symplectic basis for $V$ specifying $\Sigma$ as in Lemma
\ref{symplecticframe}, and $\widetilde{\Sigma}$ the apartment of
$\Xi_{2n}$ specified by $\cB$.
Let $D \in \widetilde{\Sigma}$ be any chamber
containing $C$. Then $D$ corresponds to the chain
$\pi L_0 \subsetneq L_1 \subsetneq \cdots \subsetneq L_n \subsetneq
L_{n + 1} \subsetneq \cdots \subsetneq L_{2n - 1} \subsetneq L_0$ for
some lattices $L_{n + 1}, \ldots, L_{2n - 1}$ in $V$.
For $0 \leq j \leq 2n - 1$, write
\[
L_j = (a_1^{(j)}, \ldots, a_n^{(j)}; b_1^{(j)}, \ldots,
b_n^{(j)})_{\cB}.
\]
\begin{nlem}
\label{chamber_special}
The two special vertices in $C$ are $[L_0]$ and $[L_n]$.
\begin{proof}
The fact that $[L_0]$ is special follows from \cite[Corollary
3.4]{shemanske} and \cite[p.\ 3411]{shemanske}.
To see that $[L_n]$ is special, note
that if $L_j$ represents a special vertex in $C$ for $1 \leq j \leq
n$, then $a_i^{(j)} + b_i^{(j)} = \mu$ for all $i$ (by \cite[Corollary
3.4]{shemanske}), where $\mu \in
\{1, 2\}$ (since $\langle L_j, L_j \rangle \subseteq \pi\cO$). But
$\mu = 2$ implies
$L_j = \pi L_0$, which is impossible.
Thus, $a_i^{(j)} + b_i^{(j)} = 1$ for all $i$ and $L_j/\pi L_0 \cong
k^n$; hence,
$j = n$.
\end{proof}
\end{nlem}

For $\cB = \{u_1, \ldots, u_n, w_1, \ldots, w_n\}$ a symplectic basis
for $V$ and $g \in \GSp_n(K)$, let
\[
\cB_g := \{\nu(g)^{-1}gu_1, \ldots, \nu(g)^{-1}gu_n, gw_1, \ldots,
gw_n\}.
\]
Note that $\cB_g$ is a symplectic basis for $V$; hence, $L = (a_1, \ldots,
a_n; b_1, \ldots, b_n)_{\cB}$ and $\ord(\nu(g)) = m$ imply $gL =
(a_1 + m, \ldots, a_n + m; b_1, \ldots, b_n)_{\cB_g}$.
\begin{nprop} \label{GSptransitive}
The group $\GSp_n(K)$ acts transitively on the
special vertices in $\Delta_n$.
\begin{proof}
Note that if $\GSp_n(K)$ acts on the special vertices in $\Delta_n$, then
\cite[Proposition 3.3]{shemanske} implies that the
action is transitive. We thus show that $\GSp_n(K)$ acts on the special
vertices in $\Delta_n$.
Let $t \in \Delta_n$ be a special vertex and $L \in t$ a representative such
that there is a primitive lattice $L_0$ with $\langle L, L \rangle \subseteq
\pi\cO$ and $\pi L_0 \subseteq L \subseteq L_0$. Let
$\Sigma$ be an apartment of $\Delta_n$ containing $t$ and $[L_0]$,
and let $\cB$ be a symplectic basis for $V$ specifying $\Sigma$ as in Lemma
\ref{symplecticframe}. Then \cite[p.\ 3411]{shemanske}, the last lemma, and
\cite[Corollary 3.4]{shemanske} imply
\[
L_0 = (c_1, \ldots, c_n; -c_1, \ldots, -c_n)_{\cB} \qquad\text{and}\qquad
L = (a_1, \ldots, a_n; \mu - a_1, \ldots, \mu - a_n)_{\cB},
\]
where $\mu \in \{1, 2\}$. Let $g \in \GSp_n(K)$ with $\ord(\nu(g)) = m$. Since
$gt = [a_1 + m, \ldots, a_n + m; \mu - a_1, \ldots, \mu - a_n]_{\cB_g}$,
\cite[Corollary 3.4]{shemanske} implies that it suffices to show $gt$ is a
vertex in $\Delta_n$. First suppose $m \equiv 0 \mod 2$, say $m = 2r$. Then
$\pi^{-r}gL_0$ is primitive, $\langle \pi^{-r}gL, \pi^{-r}gL \rangle \subseteq
\pi\cO$, and $\pi^{-r}g(\pi L_0) \subseteq \pi^{-r}gL \subseteq \pi^{-r}gL_0$;
i.e., $gt$ is a vertex in $\Delta_n$. Now suppose $m = 2r + 1$. If $\mu = 1$,
then $\pi^{-r - 1}gL$ is primitive and $gt$ is a vertex in $\Delta_n$.
Otherwise, $\mu = 2$, and $\langle \pi^{-r - 1}gL, \pi^{-r - 1}gL \rangle
\subseteq \pi\cO$. Let $\pi M_0 = (a_1 + r, \ldots, a_n + r; \mu - a_1 - r,
\ldots, \mu - a_n - r)_{\cB_g}$. Then $M_0$ is primitive and $\pi M_0 \subseteq
\pi^{-r - 1}gL \subseteq M_0$; i.e., $gt$ is a vertex in $\Delta_n$. Thus,
$\GSp_n(K)$ acts on the special vertices in $\Delta_n$.
\end{proof}
\end{nprop}

Note that by
Propositions \ref{primitivetype} and \ref{specialtype},
$[L_n]$ has type $n$. Then by Proposition \ref{Spntypes}, the type of
$[L_j]$ is in $\{n + 1, \ldots, 2n - 1\}$ for all $1 \leq j \leq n -
1$ and the type of $[L_i]$ is in $\{1, \ldots, n - 1\}$ for all $n +
1 \leq i \leq 2n - 1$.
\begin{nlem} \label{Deltachamber2}
Let $g \in \GSp_n(K)$ with $\ord(\nu(g)) \equiv 1 \mod 2$. If $L_0,
L_n, \ldots, L_{2n - 1}$ are lattices in $V$ as above, then the
vertices $g[L_n], \ldots, g[L_{2n - 1}], g[L_0]$ in $\Xi_{2n}$ are the
vertices in a chamber in $\Delta_n$.
\begin{proof}
Write $\ord(\nu(g)) = 2r + 1$. Then Lemma
\ref{chamber_special} and \cite[p.\ 3411]{shemanske} imply that
$L_n' = \pi^{-(r + 1)}gL_n$
is
primitive (see the proof of Proposition \ref{GSptransitive}).
Furthermore,
if
$L_j' = \pi^{-r}gL_j$
for $j = 0, n + 1, \ldots, 2n - 1$, then $\pi L_n' \subsetneq L_{n +
1}' \subsetneq \cdots \subsetneq L_{2n - 1}' \subsetneq L_0'
\subsetneq L_n'$ and $\langle L_j', L_j' \rangle \subseteq \pi\cO$ for
$j = 0, n + 1, \ldots, 2n - 1$; i.e.,
$[L_n'], \ldots, [L_{2n - 1}'], [L_0']$ are the vertices in a chamber
in $\Delta_n$.
\end{proof}
\end{nlem}

\begin{nlem} \label{gallery1_apt}
Let $\Sigma$ be an apartment of $\Delta_n$ and $\widetilde{\Sigma}$
the apartment of $\Xi_{2n}$ such that $\cB$ a symplectic basis for $V$
specifying $\Sigma$ implies $\cB$ specifies $\widetilde{\Sigma}$. If
$C, C'$ is a gallery in $\Sigma$, then there is a gallery $D, D'$ in
$\widetilde{\Sigma}$ such that $D$ (resp., $D'$) contains $C$ (resp.,
$C'$) and $C \neq C'$ implies $D \neq D'$.
\begin{rem}
More generally, if $C_0, \ldots, C_m$ is a gallery in $\Delta_n$, then
there is a gallery $D_0, \ldots, D_\ell$ in $\Xi_{2n}$ and integers $0
\leq i_0 < \cdots < i_m \leq \ell$ such that $D_j$ contains $C_0$ for
all $0 \leq j \leq i_0$ and $D_j$ contains $C_r$ for all $i_{r - 1} <
j \leq i_r$ and all $1 \leq r \leq m$.
\end{rem}
\begin{proof}
If $C = C'$, set
$D = D'$, where $D \in \widetilde{\Sigma}$ is a chamber containing
$C$.
Now suppose $C \neq C'$, with $C$ corresponding to the chain
\begin{equation} \label{chamber_Spn}
\pi L_0 \subsetneq L_1 \subsetneq \cdots \subsetneq L_n \subsetneq
L_0.
\end{equation}
Let $\cB$ be a symplectic basis for $V$ specifying $\Sigma$ as in
Lemma \ref{symplecticframe}, and let
$0 \leq j \leq n$ such that $C \cap C'$
corresponds to \eqref{chamber_Spn} with
$L_j$ deleted if $1 \leq j \leq n$ or with both $\pi L_0$ and $L_0$
deleted if $j = 0$. Note that
if $t'$ is the vertex in $C'$ not in $C$, then
$t'$ has a representative $L'$ such that $C'$ corresponds to
\eqref{chamber_Spn} with $L_j$ replaced by $L'$.

If
$1 \leq j \leq n  - 1$, then
\cite[p.\ 3411]{shemanske}, 
Lemma
\ref{chamber_special}, and \eqref{chamber_Spn} imply $L_0 = (a_1,
\ldots, a_n; -a_1, \ldots, -a_n)_{\cB}$ and
$L_n = (b_1, \ldots, b_n; 1 - b_1, \ldots, 1 - b_n)_{\cB}$,
where $a_i + 1 \geq b_i \geq a_i$ for all
$i$.
For $1 \leq i \leq n$, let $a_{n + i} = -a_i$ and $b_{n + i} = 1 -
b_i$.
Let $\{i_1, \ldots, i_n\}$ be the
$n$ values of $i$ such that $b_i = a_i + 1$, and for
$1 \leq r \leq n - 1$, set $L_{n + r} = (c_1, \ldots, c_n; c_{n + 1},
\ldots, c_{2n})_{\cB}$, where $c_\ell = b_\ell - 1 = a_\ell$ if $\ell
\in \{i_1, \ldots, i_r\}$ and $c_\ell = b_\ell$ otherwise. Then $L_n
\subsetneq L_{n + 1} \subsetneq \cdots \subsetneq L_{2n - 1}
\subsetneq L_0$, and letting
$D \in \widetilde{\Sigma}$ (resp., $D' \in \widetilde{\Sigma}$) 
be the simplex
with vertices the vertices in $C$ (resp., the vertices in $C'$),
together with
$[L_{n + 1}], \ldots, [L_{2n - 1}]$ finishes the proof in this case

If
$j = n$,
write $L_0 = (a_1, \ldots, a_n; -a_1, \ldots, -a_n)_{\cB}$,
$L_n = (b_1, \ldots, b_n; 1 - b_1, \ldots, 1 - b_n)_{\cB}$,
and $L' = (b_1', \ldots, b_n'; 1 - b_1', \ldots, 1 - b_n')_{\cB}$.
Note that
$a_i + 1 \geq b_i, b_i' \geq a_i$ for all
$i$ and
$b_i \neq b_i'$ for at least one value of $i$. Let $L_{n + 1} = (c_1,
\ldots, c_n; c_{n + 1}, \ldots, c_{2n})_{\cB}$, where $c_i =
\min\{b_i, b_i'\}$ and $c_{n + i} = \min\{1 - b_i, 1 - b_i'\}$ for
$1 \leq i \leq n$. Then $L_{n + 1} = L_n + L'$, so $L_n, L' \subsetneq
L_{n + 1}$ and
$[L_{n + 1} : L_n] = q = [L_{n + 1} : L']$.
An obvious modification of the second half of the last paragraph
finishes the proof in this case.

Finally, if
$j = 0$,
write $L_0 = (a_1, \ldots, a_n; -a_1, \ldots, -a_n)_{\cB}$,
$L' = (a_1', \ldots, a_n'; - a_1', \ldots, -a_n')_{\cB}$,
and $L_n = (b_1, \ldots, b_n; 1 - b_1, \ldots, 1 - b_n)_{\cB}$. Note
that
$a_i + 1, a_i' + 1 \geq b_i \geq a_i, a_i'$
for all
$i$ and
$a_i \neq a_i'$ for at least one value of $i$. Let $L_{2n - 1} = (c_1,
\ldots, c_n; c_{n + 1}, \ldots, c_{2n})_{\cB}$, where $c_i =
\max\{a_i, a_i'\}$ and $c_{n + i} = \max\{-a_i, -a_i'\}$ for
$1 \leq i \leq n$. Then $L_{2n - 1} = L_0 \cap L'$, so $L_{2n - 1}
\subsetneq L_0, L'$ and
$[L_0 : L_{2n - 1}] = q = [L' : L_{2n - 1}]$.
An obvious modification of the second half of the first paragraph
finishes the proof in this case.
\end{proof}
\end{nlem}

It will turn out to
be convenient to first prove results
about the type $0$ vertices in $\Delta_n$ and to then use the
transitive action of $\GSp_n(K)$ on the special vertices in $\Delta_n$
(see Proposition \ref{GSptransitive})
to deduce
the same results about the type $n$ vertices in $\Delta_n$.
For $g \in \GL_{2n}(K)$ and a chamber $C \in \Xi_{2n}$, abuse notation
and write $gC$
for the image of the vertices in $C$ under the action of $g$.
\begin{nprop} \label{GLngallery}
The group $\GL_{2n}(K)$ (resp., $\GSp_n(K)$) maps a gallery in
$\Xi_{2n}$ of length $m$ to a gallery in $\Xi_{2n}$ of length $m$. In
particular, if $C \neq C'$ are adjacent chambers in $\Xi_{2n}$ and $g
\in \GL_{2n}(K)$ (resp., $g \in \GSp_n(K)$), then $gC \neq gC'$ are
adjacent chambers in $\Xi_{2n}$.
\begin{proof}
Let $C_0, \ldots, C_m$ be a gallery in $\Xi_{2n}$, and let $g \in
\GL_{2n}(K)$. If $m = 0$ and
$C_0$ corresponds to the chain
$\pi L_0 \subsetneq L_1 \subsetneq \cdots \subsetneq L_{2n - 1}
\subsetneq L_0$, then $g(\pi L_0) \subsetneq gL_1 \subsetneq \cdots
\subsetneq gL_{2n - 1} \subsetneq gL_0$; i.e., $gC_0$ is a chamber in 
$\Xi_{2n}$.
If $m = 1$ and $C_0 = C_1$, then
$gC_0, gC_1$ is a gallery in $\Xi_{2n}$, so
suppose $C_0 \neq C_1$. Let
$t_0, \ldots, t_{2n - 1}$ (resp., $x_0, \ldots, x_{2n - 1}$) be the
vertices in $C_0$ (resp., in $C_1$), and let
$0 \leq j \leq 2n - 1$ such that $t_j \neq x_j$.
For $0 \leq i \leq 2n - 1$, let $L_i \in t_i$ (resp., let
$M_i \in x_i$) such that $\pi L_0 \subsetneq L_1 \subsetneq \cdots
\subsetneq L_{2n - 1} \subsetneq L_0$ (resp.,
$\pi M_0 \subsetneq M_1 \subsetneq \cdots \subsetneq M_{2n - 1}
\subsetneq M_0$) corresponds
to $C_0$ (resp., to
$C_1$).
Then
$g(\pi L_0) \subsetneq gL_1 \subsetneq \cdots \subsetneq gL_{2n - 1}
\subsetneq gL_0$ (resp.,
$g(\pi M_0) \subsetneq gM_1 \subsetneq \cdots \subsetneq gM_{2n - 1}
\subsetneq gM_0$). Since
$t_i = x_i$ implies $gt_i = gx_i$,
$gC_0, gC_1$ is a gallery in $\Xi_{2n}$. The fact that
$gC_0 \neq gC_1$ follows from the fact that $gx_j \neq gt_j$.
The proof for
$m \geq 2$ follows from the fact
that $gC_i, gC_{i + 1}$ is a gallery in $\Xi_{2n}$ for all $0 \leq i
\leq m - 1$.
\end{proof}
\end{nprop}

\subsection{Counting close vertices in $\Delta_n$}

Let
$\Gamma = \Sp_n(\cO)$, and note that
%
%
the analogues of the results in section 4.1 of \cite{shemanske} hold
if $\GSp_n(K)$ acts on the lattices in $V$ on the left (rather than on
the right).
The following is
an analogue of Theorem 3.3 of \cite{ss} for the special
vertices in $\Delta_n$.
\begin{nthm}
\label{3.3}
If $t \in \Delta_n$ is a special
vertex, then the number of vertices in $\Delta_n$ close to $t$ is the
number of left cosets of $\Gamma$ in
\[
\Gamma \diag(1, \underbrace{\pi, \ldots, \pi}_{n - 1}, \pi^2, \pi,
\ldots, \pi) \Gamma.
\]
\begin{proof}
First note that by Proposition \ref{specialtype}, a special vertex in
$\Delta_n$ has type either $0$ or $n$. Let $t \in \Delta_n$ be a
special
vertex and $t' \in \Delta_n$ a vertex close to $t$. Then there are
adjacent chambers $C, C' \in \Delta_n$ such that $t \in C$, $t' \in
C'$, but $t, t' \not\in C \cap C'$. Let $\Sigma$ be an apartment of
$\Delta_n$ containing $C$ and $C'$. If $t$ has type $0$, then by
Lemma \ref{primitiverep}, we may assume that relative to some
symplectic basis $\cB$ for $V$ specifying $\Sigma$,
$t = [0, \ldots, 0; 0, \ldots, 0]_{\cB} \in C_0$, where $C_0 \in
\Sigma$ is the chamber with
vertices $[0, \ldots, 0; 0, \ldots, 0]_{\cB}, [0, 1, \ldots, 1; 1,
\ldots, 1]_{\cB}, \ldots, [0, \ldots, 0; 1, \ldots, 1]_{\cB}$. A
straightforward modification of the fourth and fifth paragraphs of the
proof of \cite[Theorem 3.3]{ss} using the reflections defined in
\cite[p.\ 3411]{shemanske} finishes the proof in this case.
Now suppose $t$ has type
$n$,
and let
$\cB$ be a symplectic basis for $V$ specifying $\Sigma$ as in Lemma
\ref{symplecticframe}.
Let $\widetilde{\Sigma}$ be the apartment of $\Xi_{2n}$ specified by
$\cB$, and let
$D, D' \in \widetilde{\Sigma}$ be adjacent chambers with $C$ in $D$,
$C'$ in $D'$, and $D \neq D'$ as in Lemma
\ref{gallery1_apt}. Let
$g \in \GSp_n(K)$ with $\ord(\nu(g)) \equiv 1 \mod 2$. Then by Proposition
\ref{action_type},
$gt$ has type $0$. By Lemma \ref{Deltachamber2}, $gD$ (resp., $gD'$)
contains a chamber $C_1 \in \Delta_n$ (resp., a chamber $C_1' \in
\Delta_n$) with $gt \in C_1$ (resp., with $gt' \in C_1'$).
Furthermore,
$gD \neq gD'$ are adjacent chambers in $\Xi_{2n}$ and
$gt, gt' \not\in gD \cap gD'$ by Proposition \ref{GLngallery};
i.e., $gt$ and $gt'$ are close vertices in $\Delta_n$. Finally, if
%
$S_t$
and $S_{gt}$ are the sets
of vertices in $\Delta_n$ close to $t$ and $gt$, respectively, then
$\Card(S_t) = \Card(S_{gt})$, and the
last paragraph
finishes the proof.
\end{proof}
\end{nthm}
\begin{rem}
The analogues of the results in \cite[Section 4.1]{shemanske} also hold
if $\Sp_n(\cO)$ and $\GSp_n^S(K)$ are replaced by $\GSp_n(\cO) =
\GL_{2n}(\cO) \cap \GSp_n(K)$ and $\GSp_n(K)$, respectively, and with
$\GSp_n(K)$ acting on the left rather than on the right. In addition, the
analogue of the above theorem holds with $\Gamma = \GSp_n(\cO)$; hence,
so does Corollary \ref{cosetnumber}.
\end{rem}

\begin{figure}
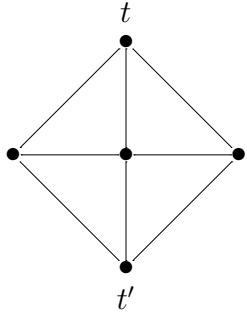

\[
\xy
(0, 15)*{\bullet}="t";
(0, 19)*{\text{$t$}};
(0, -15)*{\bullet}="b";
(0, -19)*{\text{$t'$}};
(-15, 0)*{\bullet}="l";
(0, 0)*{\bullet}="c";
(15, 0)*{\bullet}="r";
{\ar@{-} "t";"c"; };
{\ar@{-} "t";"l"; };
{\ar@{-} "t";"r"; };
{\ar@{-} "l";"c"; };
{\ar@{-} "c";"r"; };
{\ar@{-} "b";"l"; };
{\ar@{-} "b";"r"; };
{\ar@{-} "b";"c"; };
\endxy
\]
\caption{Two close special vertices, both of type $0$, in $\Delta_2$.}
\label{closeSp2}
\end{figure}

We now count the number of vertices in $\Delta_n$ close to a given
special vertex $t \in
\Delta_n$. By Proposition \ref{specialtype} and
Theorem \ref{3.3},
it suffices to assume $t$
has type $0$.
%
By Proposition \ref{primitivetype}, $t$ has a primitive representative
$L$, so a chamber $C \in \Delta_n$ containing $t$ corresponds to a
chain of the form
\begin{equation} \label{Spnchamberchain}
\pi L \stackrel{q}{\subsetneq} L_1 \stackrel{q}{\subsetneq} \cdots
\stackrel{q}{\subsetneq} L_n \stackrel{q^n}{\subsetneq} L.
\end{equation}
The codimension-one face in $C$ not containing $t$
thus corresponds to the chain
\[
L_1 \stackrel{q}{\subsetneq} \cdots \stackrel{q}{\subsetneq} L_n,
\]
and a vertex in $\Delta_n$ is close to $t$ if it has a primitive
representative $M \neq L$ such that
\begin{equation} \label{Spnclosechain}
\pi M \stackrel{q}{\subsetneq} L_1 \stackrel{q}{\subsetneq} \cdots
\stackrel{q}{\subsetneq} L_n \stackrel{q^n}{\subsetneq} M.
\end{equation}
Given the lattice $L_1$, the possible $L$ and $M$
satisfy $L \neq M \subsetneq \pi^{-1}L_1$ with $[\pi^{-1}L_1 : L] =
q =
[\pi^{-1}L_1 : M]$
and both $L$ and $M$ primitive. On the other hand,
if $t, t' \in \Delta_n$ are
close type $0$ vertices,
then there must be primitive representatives $L \in t$ and $M \in t'$
and lattices $L_1, \ldots, L_n$ as in \eqref{Spnchamberchain} such that
$L \neq M \subsetneq \pi^{-1}L_1$. The same argument as in Section
\ref{SLn_close}
shows that $\pi^{-1}L_1 = L + M$, but we can vary $L_2, \ldots, L_n$
as long as $\langle L_i, L_i \rangle \subseteq \pi\cO$ for all $2 \leq
i \leq n$ and the chains
\[
\pi L \subsetneq L_1 \subsetneq L_2 \subsetneq \cdots \subsetneq L_n
\subsetneq L \qquad\text{and}\qquad \pi M \subsetneq L_1 \subsetneq
L_2 \subsetneq \cdots \subsetneq L_n \subsetneq M
\]
correspond to chambers in $\Delta_n$. In other words (as in the case
of $\Xi_n$), if $t, t' \in \Delta_n$ are close type $0$ vertices,
there may be more than one pair of adjacent chambers $C, C' \in
\Delta_n$ such that $t \in C$, $t' \in C'$, and $t, t' \not\in C \cap
C'$ (see Figure \ref{closeSp2}).
We return to this later.

Before we count the number of vertices in $\Delta_n$ close to
$t$,
we make a few observations similar to those preceding Proposition
\ref{SLnd1count}. Fix a primitive representative $L \in t$. Then
$L/\pi L \cong k^{2n}$
is endowed with a
non-degenerate, alternating $k$-bilinear form. Moreover,
the Correspondence Theorem, the fact that any
$\cO$-submodule
of $L$ containing $\pi L$
is a lattice in $V$, and the fact that
every $1$-dimensional $k$-subspace of $L/\pi L$ is \ti\ imply that
the number of $L_1$ is
the number of $1$-dimensional $k$-subspaces of $L/\pi L$. Given $L_1$,
let $C \in \Delta_n$ be a chamber containing $[L_1]$ and $t$, and let
$A$ be the codimension-one face in $C$ not containing $t$. Then
the number of primitive lattices $M \neq L$ in $V$ such that $M
\subsetneq \pi^{-1}L_1$ and $[\pi^{-1}L_1 : M] = q$ is one less than
the number of chambers in $\Delta_n$ containing $A$.
\begin{nprop} \label{Spnd1count}
If $t \in \Delta_n$ is a special vertex, then
the number $\omega(\Delta_n)$ of vertices in $\Delta_n$ close to $t$
is
\[
\frac{q^{2n} - 1}{q - 1} \cdot q
\]
(independent of $t$).
\begin{proof}
This follows from the preceding
comments, the fact that the number of $1$-dimensional subspaces of
$\F_q^m$ is exactly $(q^m - 1)/(q - 1)$, and
Proposition \ref{codimno}.
\end{proof}
\end{nprop}
\begin{ncor} \label{cosetnumber}
The number of left cosets of $\Gamma = \Sp_n(\cO)$ in
\[
\Gamma\diag(1, \underbrace{\pi, \ldots, \pi}_{n - 1}, \pi^2, \pi,
\ldots, \pi)\Gamma
\]
is $((q^{2n} - 1) \cdot q)/(q - 1)$.
\begin{proof}
This follows from Theorem \ref{3.3} and the last proposition.
\end{proof}
\end{ncor}

Proposition \ref{chambercount} and
the last proposition prove
the following analogue of
Theorem \ref{SLnrelation}.
\begin{nthm} \label{Spnrelation}
Let $r(\Delta_n)$ be the number of chambers in $\Delta_n$ containing a
given special vertex (as in Proposition \ref{chambercount}) and
$\omega(\Delta_n)$ the number of vertices in $\Delta_n$ close to a
given special vertex in $\Delta_n$ (as in Proposition
\ref{Spnd1count}). Then for all $n \geq 2$, $q \cdot r(\Delta_n) =
r(\Delta_{n - 1})\ \omega(\Delta_n)$, where $r(\Delta_1) = q + 1$.
\end{nthm}

When the given vertex in $\Delta_n$ has type $0$, we can also give a
combinatorial proof of
Theorem \ref{Spnrelation}.
%
As
in Section \ref{SLn_close}, if $t \in \Delta_n$ is a fixed type $0$
vertex, then we can try to count the number of vertices in $\Delta_n$
close to $t$ by counting the number of galleries (in $\Delta_n$) of
length $1$ starting at a chamber containing $t$ and ending at a
chamber not containing $t$. An argument analogous to that in Section
\ref{SLn_close} shows that if $t' \in \Delta_n$ is a vertex close to
$t$, then
$\omega(\Delta_n) = (r(\Delta_n) \cdot q)/m(\Delta_n, t, t')$, where
$m(\Delta_n, t, t')$ is the number of galleries of length $1$ in
$\Delta_n$ whose initial chamber contains $t$ and whose ending chamber
contains $t'$.

To determine
$m(\Delta_n, t, t')$,
fix the following notation for the rest of this section. For close
special vertices
$t, t' \in \Delta_n$
with $t$ of type $0$,
let $L \in t$, $M \in t'$ be primitive representatives
(by Proposition \ref{primitivetype})
such that there are lattices $L_1, \ldots, L_n$ as in
\eqref{Spnchamberchain} and \eqref{Spnclosechain} with 
$\langle L_i, L_i \rangle \subseteq \pi\cO$ for all $1 \leq i \leq n$.
Recall that
$L_1 = \pi(L + M)$, but we can vary $L_2, \ldots, L_n$ as long as
$\langle L_i, L_i \rangle \subseteq \pi\cO$ for all $2 \leq i \leq n$
and the chains
\[
\pi L \subsetneq L_1 \subsetneq L_2 \subsetneq \cdots \subsetneq L_n
\subsetneq L \qquad\text{and}\qquad \pi M \subsetneq L_1 \subsetneq
L_2 \subsetneq \cdots \subsetneq L_n \subsetneq M
\]
correspond to chambers in $\Delta_n$. As in Section \ref{SLn_close},
each
gallery in $\Delta_n$ counted by $m(\Delta_n, t, t')$ is uniquely
determined by $L_2, \ldots, L_n$.
Define two vertices in $\Delta_n$ to be
{\it adjacent} if they are distinct and incident.
\begin{nlem} \label{Spnadjacent}
Let $t, t' \in \Delta_n$ be adjacent vertices such that $t$ has a
primitive representative $L$. Then $t'$ has a unique representative
$L'$ such that $\langle L', L' \rangle \subseteq \pi\cO$ and $\pi L
\subsetneq L' \subsetneq L$.
\begin{proof}
Since
$t$ and
$t'$ are adjacent vertices in $\Xi_{2n}$,
by
Proposition \ref{adjacentreps}, $t'$ has a unique representative $L'$
such that $\pi L \subsetneq L' \subsetneq L$. It thus
suffices to show that $\langle L', L' \rangle \subseteq \pi\cO$.
But $t$ and $t'$
incident vertices in $\Delta_n$ with $t \neq t'$ implies
they have representatives $M \in t$
and $M' \in t'$ such that
there is a primitive lattice $L_0$ with $\langle M, M \rangle
\subseteq \pi\cO$, $\langle M', M' \rangle \subseteq \pi\cO$, and
either $\pi L_0 \subseteq M \subsetneq M' \subseteq L_0$ or $\pi L_0
\subseteq M' \subsetneq M \subseteq L_0$. Suppose
$\pi L_0 \subseteq M \subsetneq M' \subseteq L_0$ (resp., $\pi L_0
\subseteq M' \subsetneq M \subseteq L_0$). Then $M$ and $\pi L$
(resp., $M$ and $L$) homothetic implies $\pi L = \pi^rM$ (resp., $L =
\pi^r M$) for some $r \in \Z$;
hence, $\pi L
\subsetneq \pi^rM' \subsetneq
L$. Let $L' = \pi^rM'$. Since $L$ is primitive,
$\langle \pi^{r - 1}M, \pi^{r - 1}M \rangle \subseteq \cO$ (resp.,
$\langle \pi^rM, \pi^r M \rangle \subseteq \cO$). On the other hand,
$\langle \pi^{r - 1}M, \pi^{r - 1}M \rangle \subseteq \pi^{2(r - 1) +
1}\cO$ (resp., $\langle \pi^rM, \pi^rM \rangle \subseteq \pi^{2r +
1}\cO$), so
$r \in \Z^+$ (resp., $r \in \Z^{\geq 0}$) and
$\langle L', L' \rangle
\subseteq \pi\cO$.
%
\end{proof}
\end{nlem}

Consider the set of vertices in $\Delta_n$ that are adjacent to $t,
t'$, and $[L + M]$, and define
two such vertices to be incident
if they are incident
as vertices in $\Delta_n$.
Let $\Delta_n^c(t, t')$ be the set consisting of
\begin{itemize}
\item the empty set,
\item all vertices in $\Delta_n$ adjacent to $t, t'$, and $[L + M]$,
and
\item all finite sets $A$ of vertices in $\Delta_n$ adjacent to $t,
t'$, and $[L + M]$ such that any two vertices in $A$ are adjacent.
\end{itemize}
Then $\Delta_n^c(t, t')$ is a simplicial complex.
In particular,
$\Delta_n^c(t, t')$
is a subcomplex of $\Delta_n$.
%
\begin{nlem}
If $\emptyset \neq A
\in \Delta_n^c(t, t')$ is an $i$-simplex, then $A$ corresponds to a
chain of lattices $M_1 \subsetneq \cdots \subsetneq M_{i + 1}$, where
$\langle M_j, M_j \rangle \subseteq \pi\cO$ for all $1 \leq j \leq i +
1$ and $\pi(L + M) \subsetneq M_1 \subsetneq \cdots \subsetneq M_{i +
1} \subsetneq L \cap M$. In particular, $A$ has at most $n - 1$
vertices.
\begin{proof}
As in the proof of Lemma \ref{SLnclosesimplex}, we proceed by
induction on $i$. If $i = 0$, then
$L$
primitive,
$A$
adjacent to $t$, and Lemma \ref{Spnadjacent} imply $A$
has a unique representative $M_1$
such that $\langle M_1, M_1
\rangle \subseteq \pi\cO$ and $\pi L \subsetneq M_1
\subsetneq L$. Since $A$
and $[L + M]$ are adjacent vertices in $\Xi_{2n}$,
either $M_1
\subsetneq \pi(L + M)$ or $M_1
\supsetneq \pi(L + M)$ by \cite[p.\ 322]{garrett}. But $M_1
\subsetneq \pi(L + M)$ means
$\pi L \subsetneq M_1
\subsetneq \pi(L + M)$, which is impossible since $[\pi(L + M) : \pi
L] = q$; hence, $M_1
\supsetneq \pi(L + M)$. Then $A$
and $t'$
adjacent vertices in $\Xi_{2n}$ and
\cite[p.\ 322]{garrett} imply
that either $M_1
\subsetneq M$ or $M_1
\supsetneq M$. Since $M_1
\supsetneq M$ means
$M \subsetneq M_1
\subsetneq L$, which contradicts the fact that
$[M : \pi(L + M)]
= [L : \pi(L + M)]$, $M_1
\subsetneq M$ and
$M_1 \subseteq L \cap M$. Moreover, $\langle M_1, M_1
\rangle \subseteq \pi\cO$ implies $M_1/\pi L$
is a \ti\ $k$-subspace of $L/\pi L$ and $[M_1
: \pi L] \leq q^n$. The fact that
$[L \cap M : \pi L] = q^{2n - 1}$ finishes the proof in this case.

Recall
that $\langle \cdot, \cdot \rangle$ induces
a non-degenerate, alternating $k$-bilinear form
on $L/\pi L$. Then
with respect to this induced bilinear form,
$(L \cap M)/\pi L$ is the orthogonal complement of $\pi(L + M)/\pi L$
in $L/\pi L$.
In addition,
$\langle \cdot, \cdot \rangle$ induces a non-degenerate, alternating
$k$-bilinear form
on
$(L \cap M)/\pi(L + M) \cong k^{2(n - 1)}$, and
there
is a bijection between nested sequences $S_1 \subsetneq \cdots
\subsetneq S_{i + 1}$ of \ti\ $k$-subspaces of $(L \cap M)/\pi(L + M)$
and chains of $\cO$-submodules $M_1 \subsetneq \cdots \subsetneq M_{i
+ 1}$ of $L \cap M$ containing $\pi(L + M)$ with $\langle M_j, M_j
\rangle \subseteq \pi\cO$ for all $1 \leq j \leq i + 1$.
An obvious modification of the second paragraph of the proof of Lemma
\ref{SLnclosesimplex} finishes the proof.
\end{proof}
\end{nlem}

Recall that $\Delta_n^s(k)$ denotes the spherical $C_n(k)$ building
described in \cite[pp.\ 5 -- 6]{ronan}.
\begin{nprop} \label{Spnclosebldg}
For any close special vertices $t, t' \in \Delta_n$ with $t$ of type
$0$, $\Delta_n^c(t, t')$
is isomorphic (as a poset) to
$\Delta_{n - 1}^s(k)$ (independent of $t$ and $t'$ with $t$ of type
$0$).
\begin{proof}
Let $L \in t, M \in t'$ be primitive representatives as in the
paragraph preceding Lemma \ref{Spnadjacent}, and let
$\Delta_{n - 1}^s(k)$ be the spherical $C_{n - 1}(k)$ building with
simplices the empty set, together with
the nested sequences
of non-trivial, \ti\ $k$-subspaces of $(L \cap M)/\pi(L + M)$. Then
the last lemma implies that
there is a bijection between the $i$-simplices in $\Delta_n^c(t, t')$
and the $i$-simplices in $\Delta_{n - 1}^s(k)$ for all $i$. Since
this bijection preserves the partial order (face) relation, it
is a poset isomorphism.
\end{proof}
\end{nprop}
\begin{nprop} \label{Spnrelation2}
If $t, t' \in \Delta_n$ are close special vertices with $t$ of type
$0$, then
$m(\Delta_n, t, t') =
r(\Delta_{n - 1})$ (independent of $t$ and $t'$). In particular,
$\omega(\Delta_n) = (r(\Delta_n) \cdot q)/r(\Delta_{n - 1})$.
\begin{proof}
The proof is an obvious modification of the proof of Theorem
\ref{SLnrelation2}.
\end{proof}
\end{nprop}

%
\bibliographystyle{abbrv}
\bibliography{blah}

\begin{thebibliography}{1}

\bibitem{brown}
K.~Brown.
\newblock {\em Buildings}.
\newblock Springer-Verlag, New York, 1989.

\bibitem{cartwright}
D.~Cartwright.
\newblock Harmonic functions on buildings of type $\widetilde{A}_n$.
\newblock In M.~Picardello and W.~Woess, editors, {\em Random Walks and
  Discrete Potential Theory}, pages 104 -- 138, Cambridge, 1999. Cambridge
  University Press.

\bibitem{garrett}
P.~Garrett.
\newblock {\em Buildings and Classical Groups}.
\newblock Chapman \& Hall, London, 1997.

\bibitem{parkinsonja}
J.~Parkinson.
\newblock Buildings and {H}ecke algebras.
\newblock {\em Journal of Algebra}, 297(1):1 -- 49, 2006.

\bibitem{ronan}
M.~Ronan.
\newblock {\em Lectures on Buildings}, volume~7 of {\em Perspectives in
  Mathematics}.
\newblock Academic Press, Inc., Boston, MA, 1989.

\bibitem{ss}
A.~Schwartz and T.~Shemanske.
\newblock Maximal orders in central simple algebras and {B}ruhat-{T}its
  buildings.
\newblock {\em Journal of Number Theory}, 56(1):115 -- 138, 1996.

\bibitem{shemanske}
T.~Shemanske.
\newblock The arithmetic and combinatorics of buildings for $\text{Sp}_n$.
\newblock {\em Transactions of the American Mathematical Society}, 359(7):3409
  -- 3423, 2007.

\bibitem{tits}
J.~Tits.
\newblock Reductive groups over local fields.
\newblock In A.~Borel and W.~Casselman, editors, {\em Automorphic Forms,
  Representations and $L$-functions. Part $1$}, pages 29 -- 69, Providence, RI,
  1979. American Mathematical Society.

\end{thebibliography}
\end{document}